# UNIFIED FOUNDATIONS FOR MATHEMATICS

Mark Burgin

Department of Mathematics
University of California, Los Angeles
405 Hilgard Ave.
Los Angeles, CA 90095, USA

There are different meanings of foundation of mathematics: philosophical, logical, and mathematical. In this paper, foundations of mathematics are considered as a theory that provides means (concepts, structures, methods etc.) for the development of whole mathematics, or at least, of its major part. Set theory has been for a long time the most popular foundation. However, it was not been able to win completely over its rivals: logic as the base for logicism, theory of algorithms and computation as the base for constructive approach to mathematics, and a more recent candidate for foundations – the theory of categories. Moreover, practical applications of mathematics and its inner problems caused creation of different generalization of sets: multisets, fuzzy sets, rough sets etc. Thus, the whole situation brings us to a problem: Is it possible to find the most fundamental structure in mathematics or fundamental mathematical entities must be as diverse as mathematics itself? The situation is similar to the quest of physics for the most fundamental "brick" of nature and for a grand unified theory of nature.

This work is aimed at a demonstration that in contrast to physics, which is still in search for a unified theory and is striving to get particles that are the most elementary, in mathematics such a theory exists. It is the theory of named sets. It is demonstrated in this paper that the construction of a named set allows one to systematize all approaches to foundation of mathematics. Evidence is provided that named set is the most fundamental concept and structure of mathematics as it underlies other fundamental concepts and structures: sets, logical, algorithmic/computational, and categorical structures and entities. Although named sets are closely related to ordinary sets, existence of independent axiomatics of the theory of named sets demonstrates that named sets are not only conceptually more fundamental than sets, but also are independent in the logical setting.

**1. Introduction**

Foundations of mathematics, like foundations of any other scientific discipline, deal with three kinds of problems:

1. What objects are studied and which of them are basic for mathematics.
2. What are basic properties of the studied objects. These properties have usually the form of principles, postulates or axioms.
3. How other objects are build from basic objects and how other properties are derived from basic properties.

In this paper, only the first problem is considered because it is the most fundamental and related to the problem of knowledge unification. Unification of knowledge is one of the prior goals in mathematics and science. An important aspect of such unification is a search for the most fundamental object both for a discipline and for its object domain. Sometimes such search is easy and the fundamental structure is implanted into the theory from the beginning. As an example, we can take theory of groups. In it group is the fundamental object and everything else is based on it. In theory of categories, category is such a fundamental object. In theory of gauge fields, gauge field is a fundamental object. However, when we take a more general and extended area like physics or mathematics, the problem of finding the most fundamental object becomes rather complicated.

For example, the most fundamental object for physics has to be a model of the most fundamental object in nature if we assume that physics is the most basic natural science. We know that for thousands of years the best thinkers of mankind tried to solve the problem of finding the most fundamental object for nature. At first, it was a quest of philosophers who called such fundamental objects atoms and suggested different fantastic ideas on their properties. Then physicists began their search for the most fundamental object in nature. At first, they found molecules, then atoms, then elementary particles, and then quarks. There were many achievements along this way, but no solution has been found yet.

The situation in mathematics is, in some sense, dual to the situation in physics. While in physics the results of the search for the most fundamental element are very vague with a lot of particles and physical fields, mathematicians solved the problem of finding the most fundamental object for mathematics. However, the solution was not unique because it depended on some metamathematical assumptions.

Although the majority of mathematicians accept set theory as the unique foundations of mathematics, there are other fields that are also suggested as such foundations. The most important of them are: category theory as dynamic substitute for sets (Goldblatt, 1979); logic in logicism; and the theory of algorithms in constructivism. In addition, to overcome limitations of the set theoretical representation of reality, such approaches as topology without points (Menger, 1940) and mereology (Leonard, and Goodman, 1940; Leśniewski, 1992) were introduced into mathematics. At the same time, sets as they are considered in mathematics are challenged by some more general constructions that are based on physical considerations and practical intuition. The most popular of them are: multisets (Aigner, 1979; Knuth, 1997; 1998), fuzzy sets (Zadeh, 1965), rough sets (Pawlak, 1982; 1991) and some others. In spite of the fact that models of these systems are constructed by means of sets, intuition that gave birth to them suggests that in general these objects cannot be reduced to sets. Construction of these objects as systems of sets and mappings is due to the fact that almost all mathematicians think in terms of sets and their mappings.

It is necessary to remark that traditional usage of terminology based on terms "element" and "set" creates an illusion that set theory is the ultimate base for mathematics. The situation is similar to the illusion of ancient Greeks who thought that it is possible to reduce all mathematics to arithmetic of natural numbers. The development of mathematics demonstrated that both illusions are invalid.

Such a situation reveals two problems:

*Is it possible to find more general foundations than sets such that these foundations encompass all other approaches?*

*Is there a mathematical object that comprises all generalizations of sets?*

This paper suggests a positive answer to both questions. To clarify the situation with sets in mathematics, we take here, as a premise, evidence given by different researchers that sets represent only one aspect of mathematics and consider a more general than sets structure, *named sets*. The goal is to demonstrate that named sets form a broader foundation of mathematics and provide a context that encompasses those mathematical structures that either are not comprised by the set theoretical setting (e.g., categories) or are included into the scope of the set theoretical approach only by artificial means (e.g., multisets and fuzzy sets). In doing this, we encounter an interesting phenomenon. Contrary to the first impression that implies a possibility to build named sets using ordinary sets, it is discovered that it is not true: while ordinary sets are specific named sets, in which all elements have the same name that unites all these elements into one set. Exposition here makes emphasis on informal mathematical reasoning, assuming that formal systems of logic are built to represent less formal mathematical reality, and not vice versa. The goal is to consider mathematics in a broader context of the whole world in contrast to defining mathematics as knowledge about formal systems built of sets. In addition, we also give formal arguments to demonstrate that named sets (fundamental triads) provide the most fundamental foundation for mathematics.

In the second section, going after Introduction, named sets, which are also called fundamental triads, are introduced. There are three approaches to introduction of named sets: informal descriptive, formal constructive, and axiomatic. Here we consider only the first and the second approaches. It is possible to find an axiomatic approach in (Burgin, 1993; 2004).

The third section contains analysis of different approaches to foundations of mathematics in the context of named sets. In Subsection 3.1, relations between ordinary and named sets are studied. It is demonstrated that ordinary sets are special cases of named sets. In Subsection 3.2, relations between various generalizations of ordinary sets (such as multisets, fuzzy sets, and rough sets) and named sets are studied. It is demonstrated that all considered generalizations of ordinary sets are special cases of named sets. In Subsection 3.3, relations between structures of mathematical logic and named sets are studied. It is demonstrated that basic structures of logic are special cases of named sets. In Subsection 3.4, relations between mathematical categories and named sets are studied. It is demonstrated that categories are systems of specific named sets. In Subsection 3.5, such basic mathematical object as natural numbers is considered. It is demonstrated that natural numbers are properties of sets, while properties are specific named sets. In Subsection 3.6, relations between algorithms and named sets are studied. It is demonstrated that basic algorithmic structures are special cases of named sets.

## 2. The Concept of a Named Set

An attempt to extract the common essence of those mathematical constructions that are more general than sets (fuzzy sets, multisets, *L*-fuzzy sets, rough sets, genuine sets, etc.) brings us to the following structure:

$$\textbf{Entity 1} \xrightarrow{\textbf{connection}} \textbf{Entity 2} \qquad (1)$$

or

$$\textbf{Essence 1} \xrightarrow{\textbf{correspondence}} \textbf{Essence 2} \qquad (2)$$

This structure is called a *named set* or *fundamental triad*.

The term *named set* points out that it is possible to build models of named sets in a set theoretical context as the following diagrams demonstrate:

$$\textbf{Set 1} \xrightarrow{\textbf{correspondence}} \textbf{Set 2} \qquad (3)$$

or

$$\textbf{Collection 1} \xrightarrow{\textbf{relation}} \textbf{Collection 2} \qquad (4)$$

Historically, named sets were discovered in 1981 in an attempt of unification and merging the constructions of multisets and fuzzy sets that generate two orthogonal directions in set theory extensions. At first, named sets were considered as an abstract mathematical essence that was constructed from sets and their relations. Only later it was discovered and demonstrated that named sets exist independently from sets and are the most fundamental entities in the whole mathematics.

There are different kinds of named sets: set theoretical, categorical, mereological, etc. The most customary are set theoretical named sets because all mathematicians are accustomed to talk set theoretical language.

A visual representation of a set theoretical named set is given in Figure 1.

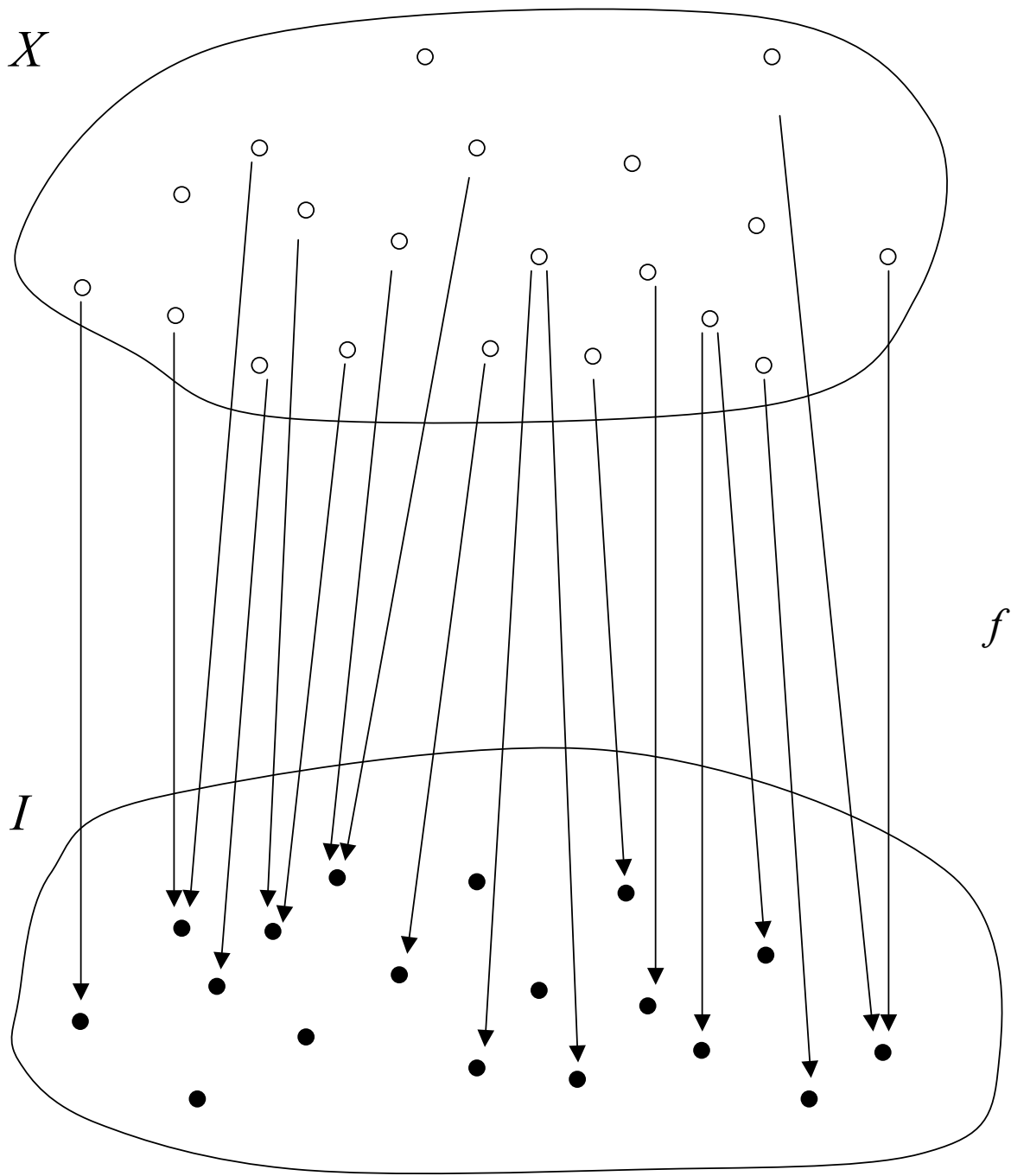

**Fig. 1.** A set theoretical named set **X** = (*X, f, I*)

Thus, a named set consists of three components and each component plays its unique role and has a specific name. In a set theoretical named set **X**, *X* and *I* are some sets, collections, or classes and *f* is a binary relation between *X* and *I*. The set (collection or class) *X* is called the *support*, the set (collection or class) *I* is called the *set of names*, and *f* is called the *naming relation* of the named set **X**.

Mereological named sets are rather different because mereology gives a very different picture of the world in comparison with the set theoretical representation. While the set theoretical approach portrays the world as a collection of sets (may be additionally structured) that consist of elements, mereology assumes that objects consist of parts. While elements are necessarily separate from one another, parts may be essentially connected, having common subparts. For instance, fingers and palms are parts of hands and hands are parts of a body. However, we cannot say that fingers and palms are elements of hands and hands are elements of a body. The membership relation that is basic for set theory is dissimilar to the very notion of parthood that mereology is about.

The roots of mereology can be traced back to the early days of philosophy, beginning with the Presocratic atomists and continuing throughout the writings of Plato, Aristotle, and Boethius. Mereology has also occupied a prominent role in the writings of medieval ontologists and scholastic philosophers such as Garland the Computist, Peter Abelard, Thomas Aquinas, Raymond Lull, and Albert of Saxony. Later this approach to reality was utilized and advanced by Leibniz and Kant. As a formal theory of parthood relations, however, mereology made its way into modern philosophy mainly through the work of Franz Brentano and of his pupils, especially Husserl. The latter may rightly be considered the first attempt at a rigorous formulation of the theory, though in a format that makes it difficult to disentangle the analysis of mereological concepts from that of other ontologically relevant notions (such as the relation of ontological dependence). It is not until Leśniewski (1916; 1992) that the pure theory of part-relations as we know it today was given an exact formulation. However, because Leśniewski's work was largely inaccessible to non-speakers of Polish, it is only with the publication of Leonard and Goodman's *The Calculus of Individuals* (1940) that this theory has become a chapter of central interest for modern ontologists and metaphysicians. Both Leśniewski's and Leonard and Goodman's original theories betray a nominalistic stand, resulting in a conception of mereology as an ontologically parsimonious alternative to set theory.

A visual representation of a mereological named set is given in Figure 2.

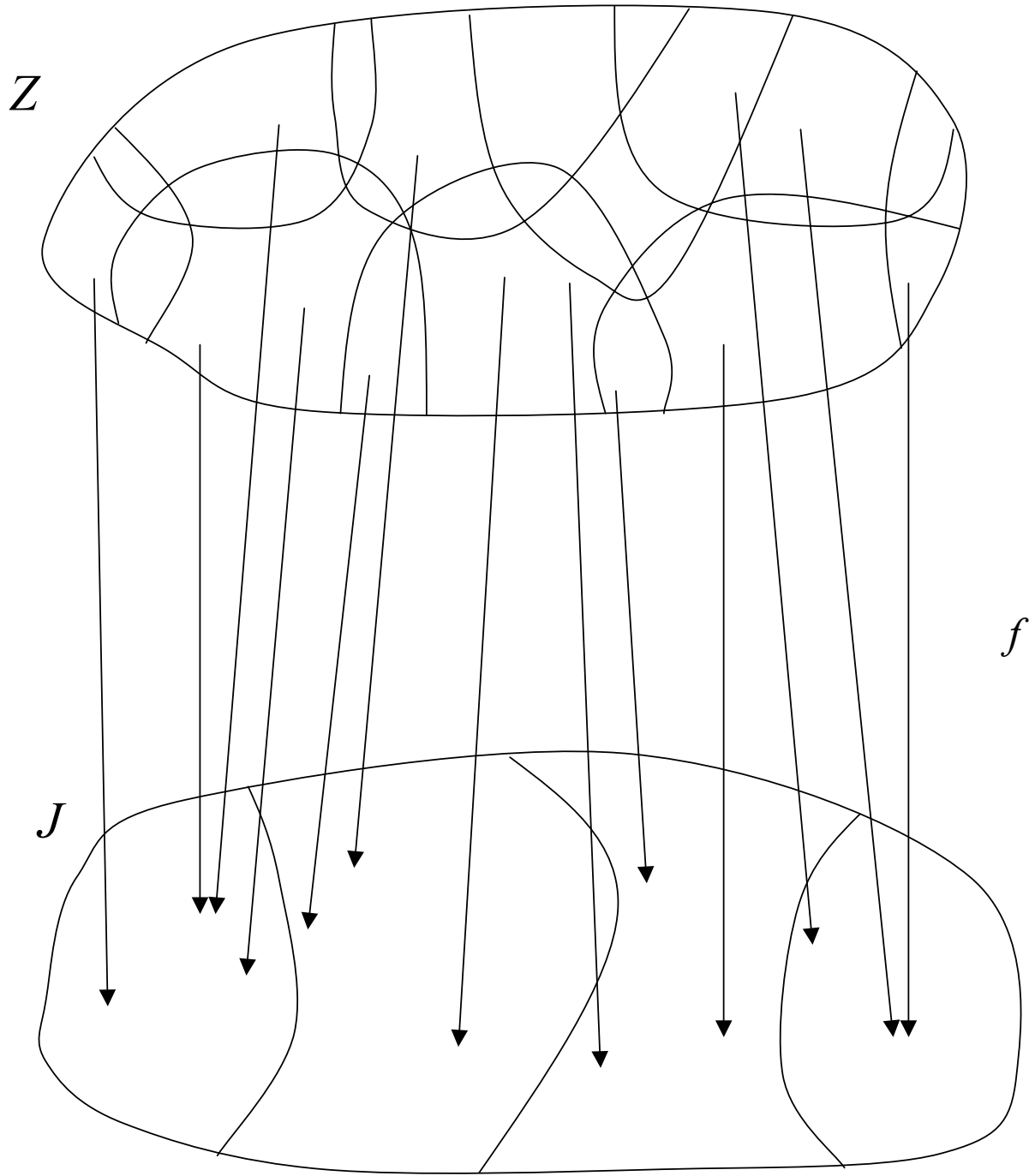

**Fig. 2.** A mereological named set **Z** = (Z, f, J)

Although the terms "named set" and fundamental triad" denote the same structure and are used interchangeably, they imply distinct intuition related to this structure, provide different perspective on it, and suggest dissimilar exposition of its properties. While the term *named set* implies a structured image, the term *fundamental triad* involves analytical representation of the structure of a named set in a form ($X, f, I$), as well as the fact that the nature of fundamental triads can be very different from the set structure. For instance, fundamental triads can be elementary indivisible objects, which consist of three parts or components. However, here "consists of" has a very different meaning from the assertion that a set consists of elements. Components of such indivisible fundamental triads do not exist outside/without those triads to which they belong.

For people, it hard to imagine how it is possible that two, three or more entities form an indivisible unity. Let us consider a visual example. Imagine a stick that is situated horizontally in front of you or take a stick and put it horizontally in front of you. Then this stick consists of:

*the right end of the stick*;

*the left end of the stick*;

and

*the part of the stick between the right end and the left end*.

We can see these three objects. We can even feel them with our fingers. So, according to principle of the contemporary physics, they exist objectively as anybody can repeat this experiment. However, it is impossible to separate these three objects. Even if somebody cuts the stick in two or more parts, each part will have the same structure. This is the structure of the fundamental triad. Thus, dividing this fundamental triad, we get other fundamental triads. This experiment gives one more proof of the fact that fundamental triads exist in reality. It is possible to find other proofs in (Burgin, 1997). In addition this experiment proves existence of indivisible, inseparable fundamental triads (named sets).

Each component of a named set (fundamental triad) plays its unique role and has a specific name. However, difference in perspective and intuition implied by these terms is supported by utilization of distinct names for the components. In a named set **X** that has the form (*X, f, I*), the following names are used to distinguish its components: *X* is called the *support*, *I* is called the *component of names*, and *f* is called the *naming correspondence* of the named set **X**. At the same time, taking representations given in Diagrams (1) and (2) and reflecting analytic role of the term "fundamental triad," we call the Entity 1 (Essence 1) the *support*, the Entity 2 (Essence 2) the *reflector*, and the connection (correspondence) the *reflection* of the fundamental triad.

People meet fundamental triads or named sets constantly in their everyday life as well as they can found them in many fundamental structures of the universe. They can even see some of fundamental triads.

Many mathematical systems are particular cases of named sets or fundamental triads. The most important of such systems are fuzzy sets (Zadeh, 1965) and multisets (Aigner, 1979; Knuth, 1997; 1998) because they are both natural generalizations of the concept of a set and have many important applications. Moreover, any ordinary set is, as a matter of fact, some named set, and namely, a singlenamed set, i.e., such a named set in which all elements have the same name (cf. Section 3.1).

It is necessary to remark that many kinds of triads are essentially indecomposable, indissoluble units. As a consequence, in such triads the objects *X, f, I* do not exist independently outside the triad which includes them.

Other kinds of triads are decomposable into three parts. But all these parts are either fundamental triads themselves or they consist of fundamental triads. It is an interesting problem whether such a decomposition may be infinite or it always ends after some (may be very big but finite) number of steps.

Named sets and fundamental triads are interchangeable as concepts. More exactly, they constitute two aspects of the same concept. Really, taking a named set **X** = (*X, f, I*), in which *X* and *Y* consist of a single element, we get a fundamental triad. At the same time, taking a fundamental triad **A** = (*A, c, B*), in which *A* and *B* are sets, collections or classes and *c* is a correspondence between their elements, we get a named set.

Thus, the name "a fundamental triad" reflects unity, wholeness, and elementarity, while the name "a named set" represents decomposability, multiplicity, and fundamentality, implying possible existence of the inner structure of the components. In what follows, we use terms "a named set" and "a fundamental triad" as synonyms related to one and the same fundamental structure. However, a fundamental triad is the name used mostly in an informal context, while a named set is usually considered as the exact mathematical structure.

At first, we give a formal mathematical definition for set theoretical named sets, utilizing category theoretical language. Let us consider three collections **Ens**, **Set**, **Col**. Each of them consists of sets (that are called objects) and their morphisms. Ob**Ens**, Ob**Set**, and Ob**Col** denote totalities of sets from **Ens**, **Set**, and **Col**, respectively. Totalities of morphisms (i.e., mappings or binary correspondences between sets) from **Ens, Set,** and **Col** are denoted by Mor**Ens**, Mor**Set**, and Mor**Col**, respectively. If $X, Y \in$ Ob **K**, then all morphisms from $X$ to $Y$ are denoted by Mor$_K$ $(X, Y)$ where K is one of the following **Ens, Set, or Col**.

Suppose that the following conditions are valid:

1) Ob**Ens**, Ob**Set** $\subseteq$ Ob**Col**;

2) Mor**Ens**, Mor**Set** $\subseteq$ Mor**Col**;

3) the totality Mor**Col** is closed in respect to the product of morphisms (it means that if $a, b \in$ Mor**Col** and their product $ab$ is defined, then $ab \in$ Mor**Col**).

Let us take some subcollection $M$ from the collection Mor**Col**. This selection may be to some extent arbitrary that makes possible, using different conditions on $M$, to define constructions necessary in each concrete case.

**Definition 2.1.** *A named set* (*in* (**Ens**, **Set**) *with respect to M*) *or an* N-*set is a triad* $\mathbf{X} = (X, r, I)$ where $X \in$ Ob**Ens**, $I \in$ Ob**Set**, $r \in$ Mor $(X, I)$ *and* $r \in M$.

$$X \xrightarrow{r} I \qquad (5)$$

**Remark 2.1.** This is a constructive or generic definition of named sets. An axiomatic definition may be found in (Burgin, 1993). While the constructive definition uses for simplicity such popular mathematical concepts as sets and their morphisms (mappings, functions, or binary relations), the axiomatic definition is completely independent from other mathematical structures to the same extent as axiomatic definitions of ordinary sets are independent from other mathematical concepts.

**Remark 2.2.** Both constructive and axiomatic definitions of named sets use concepts from other mathematical theories. Constructive approach is developed inside the theory of categories and theory of sets and classes, while the axiomatic design is, as always, based on mathematical logic. This causes an illusion that constructions of logic, sets, and/or categories are more fundamental than the structure of a named set. Existence of two types of definitions, constructive and axiomatic, demonstrates that this is not true. However, the concept of a named set has the same level of fundamentality as basic logical concepts and the concept of a category. Exact definitions of levels of fundamentality may be found in (Burgin, and Kuznetsov, 1991). However, if we consider the notion of a named set, which may be considered informally, then we will see that it is more profound than all logical and categorical notions in mathematics. The reason is that the named set structure is the most fundamental in the universe of all structures (Burgin, 1997).

Let $\mathbf{X} = (X, r, I)$ be a named set.

**Definition 2.2.** 1) *The set $X$ is called the support of* $\mathbf{X}$ *and is denoted by* $S(\mathbf{X})$; 2) *the set $I$ is called the reflector or the set of names of* $\mathbf{X}$ *and is denoted by* $N(\mathbf{X})$; 3) *the set* $N_f(\mathbf{X}) = \{ a \in I; (x \in S(\mathbf{X}) \ \& \ (x, a) \in r \}$ *is called the set of nonvoid* (*factual*) *names of the named set* $\mathbf{X}$; 4) *the mapping* (correspondence) *$r$ is called the naming mapping* (*correspondence*) *or reflection of* $\mathbf{X}$ *and is denoted by* $n(\mathbf{X})$; 5) the element $r(x) \in I$ (*for a single valued mapping $r$*) *and the set* $r(x) = \{ a \in I; (x, a) \in r \}$ (*for a binary relation or multivalued mapping $r$*) *is called the complete name of $x \in X$ in* $\mathbf{X}$; *if $r$ is not a single valued correspondence, then any element $a \in r(x)$ is called a partial name of $x$ in* $\mathbf{X}$. *Otherwise, it is called simply a name of $x$*.

**Theorem 2.1.** *If the classes* **Ens, Set**, *and* **Col** *are categories, then the collection* **NSet** *of all named sets and their morphisms is also a category*.

**Definition 2.3.** *A named set* $\mathbf{X}$ *is called*: 1) normalized *if* $N_f(\mathbf{X}) = N(\mathbf{X})$; 2) *a* singlenamed set *if* $N_f(\mathbf{X})$ *consists of a single element*; 3) *an* individually named set *if* $n(\mathbf{X})$ *is a bijection*.

**Definition 2.4.** *A named set* **Y** = (Y, b, J) *is called a named subset* (*a weak named subset*) *of the named set* **X** = (X, a, I) *if* $Y \subseteq X$, $J \subseteq I$ *and* $b = a|_{(Y,J)}$ (*and* $b \subseteq a \cap (Y \times J)$). *Such relation between named sets is called the inclusion* (*the weak inclusion*) *and is denoted by* **Y** $\subseteq$ **X** (**Y** $\subseteq_w$ **X**).

**Definition 2.5.** *If* **X** = ( X, r, I ) *and* **Y** = ( Y, q, J ) *are named sets then a morphism from* **X** *to* **Y** *is a pair* **F** = (f, g) *where* $f: X \to Y$, $g: I \to J$ *for which the equality* $fq = rg$ *is valid, i.e., the following diagram is commutative*

$$\begin{array}{ccc} I & \xrightarrow{g} & J \\ {\scriptstyle r}\uparrow & & \uparrow{\scriptstyle q} \\ X & \xrightarrow{f} & Y \end{array} \qquad (6)$$

For morphisms **F** = (f, g): X → Y and **J** = (t, s): Y → Z, their product (composition) is defined in a natural way as **FJ** = (ft, gs): X → Z.

Informally, a morphism of named sets is, as usual, a mapping or relation that preserves the naming structure.

## 3. Foundations of Mathematics and Named Sets

We have discussed in Introduction that, in contrast to physics where fundamental entities are multiplying like Fibonacci's rabbits, it was possible to discover the most basic structure of mathematics, which comprises all other foundations of mathematics in such a way that they represent different aspects of this structure. We know that centuries of search for the most fundamental object in physics resulted in discoveries of smaller and smaller particles: from molecules to atoms to elementary particles to subatomic particles, including quarks.

At the same time it has been discovered that named set plays the role of the most fundamental essence in mathematics. However, we do not mean to imply that named set is the unique fundamental structure in mathematics (and beyond, cf. Burgin, 1990; 1997), far from it. There are other fundamental structures: sets, categories, etc. But it is possible to prove that named set is the most fundamental structure in the following sense:

**The most fundamental object of mathematics is a named set because all other objects in mathematics are either some kinds of named sets or some constructions that are built from named sets**

In this paper, we argue in favor and give a variety of positive evidence for this statement. To do so, we consider different foundations of mathematics: set theory and its generalizations, mathematical logic, the theory of algorithms, theory of categories, and arithmetic of natural numbers, and show that formally or informally all of them are based on named sets. Such a unification, and not only interpretations as argues Boas (1981), helps to make mathematics more intelligible.

In particular, we explicate an interesting fact. As some of named sets are set theoretical, many have an impression that it is possible to build named sets, in general, from ordinary sets. In contrast to this illusion inspired by a set theoretical stratification of the world, we find that it is not the case that sets are more fundamental than named sets. In the domain of named sets, ordinary sets occupy a specific place, being special cases of general named sets.

Indeed, any set that is considered either in mathematics or in real life has a name. As Poincare (1908) wrote, without a name, no object exists in science or mathematics. The name of a set forms this set from separate elements by attaching to each of them a common name. Namely, all elements from a set $X$ have the common name "an element from the set $X$." As a result any set implicitly exists as a singlenamed set. In conventional set theories, supports of these singlenamed sets are considered as autonomous objects, giving the abstraction of ordinary sets. Thus, when we write that a set $X$ is equal to a set $Y$, it means that these singlenamed sets have the same support, and in many cases, it is possible to ignore existence of element names.

However, there are situations when the singlenamed structure matters. One example is from mathematics. A formal language with an alphabet $A$ is traditionally defined as an arbitrary subset of the set $A^*$ of all words in $A$. Given a Turing machine $T$, its language *L*(*T*) is often determined as the set of words accepted by *T*. However, given two Turing machines *T* and *Q*, it is impossible to determine in a general case whether *L*(*T*) = *L*(*Q*) as sets. Consequently, it is more natural and efficient to consider languages of Turing machines as singlenamed sets, which in some cases have the same support.

Another example is from biology. Interpretation of classes of biological objects not as sets but as named sets, allows one to solve the, so-called, Greg's paradox in taxonomy (Burgin, 1983).

Below we consider main directions in foundations of mathematics: set theory and its generalizations (fuzzy set theory and multiset theory), mathematical logic, the theory of algorithms, theory of categories, and arithmetic of natural numbers. For each of these directions, it is demonstrated that an attempt to extract a basic element brings us to a named set.

It is necessary to remark that some regard the theory of algorithms as identical to the theory of recursion and include it into mathematical logic. Here we distinguish mathematical logic and theory of algorithms, treating the first as a mathematical theory of reasoning/thinking and the latter as a mathematical theory of action and operation. As a result, the theory of recursion is treated as a part of theory of algorithms.

### 3.1. Set Theory as a Foundation

The notion of a *set*, *class*, collection, or aggregate is very popular in everyday life and consequently, in natural languages. In a similar way, these notions, as well as the concept of a set are among the most fundamental in mathematics. As a consequence, the most popular foundation of mathematics is axiomatic set theory.

For a long time, mathematics used only notions of a set, class, collection, or aggregate without trying to build corresponding concepts. Only in the 19$^{th}$ century the great mathematician Georg Cantor built a set theory, making set a legal concept of mathematics. Although his definition was later developed and made more strict and rigorous, it and different its versions has been used as an informal mathematical concept to our time.

Building his revolutionary set theory, Cantor (1895) gives the following definition of a set:

**Definition 3.1.** "*We understand a "set" as a combination in some whole M of definite well distinctive things m from our observation or thought. These things are called "elements" of the set M.*

Symbolically we express this as follows:

$$M = \{m\}"$$

In this definition, we encounter two named sets:

$$\textbf{a set } M \xrightarrow{\text{naming relation}} \textbf{a name "}M\textbf{" of the set } M \quad (7)$$

and

$$\textbf{elements from the set } M \xrightarrow{\text{naming relation}} \textbf{common name "}m\textbf{" of these elements} \quad (8)$$

However, other mathematicians produced and utilized a variety of informal definitions of sets. Let us consider some of them, which are used in textbooks, journal papers, and monographs in mathematics.

Courant and Robbins (1969) state that "*the concept of a class or set is one of the most important in mathematics and that a set is defined by any property or attribute* **A** *which each object considered must either posses or not posses.*"

Thus, we encounter a new named set, which is connected with the initial set:

$$\textbf{elements from a set } M \xrightarrow{\text{attributive relation}} \textbf{property or attribute A} \quad (9)$$

and a new named class:

$$\textbf{elements from a class } M \xrightarrow{\text{attributive relation}} \textbf{property or attribute A} \quad (10)$$

In some books, sets are not defined. However, they are always named. As Poincare (1908) wrote, without a name, no object exists in science or mathematics. Thus, we always have the following named set:

$$\textbf{a set} \xrightarrow{\text{correspondence}} \textbf{a name of the set} \quad (11)$$

In the informal set theory and in any other set theory (naive, constructive, descriptive, etc.), all sets are intrinsically connected with two named sets. The first kind of these named sets emerges through naming because any set has some name. As emphasized Poincare, we cannot speak about, use or construct any set without giving a name to it. Consequently, the named set (11) always exists when we have a set.

Moreover, we can see that any set is, actually, some named set as such representation more completely displays the structure of sets. Really, there are two familiar and natural ways of constructing sets: extentional and intentional (Fraenkel and Bar-Hillel, 1958). On the one hand, in the extentional approach, given a multiplicity of objects, some or all of these objects can be conceived together as forming a set. While doing so, we baptize this set by some name. The simplest name of a set has the form of a letter $X$ or $A$. In general, a name of a set is a linguistic expression $P$ or, in particular, some letter (like $A$ or $B$) in the naive set theory of Georg Cantor. A name of a set in descriptive set theory is some descriptive representation using constructive (in some sense) operations. A name of a set in constructive or recursive set theory is an algorithm generating (computing) the elements from a constructive set or an identification algorithm of the elements from the basic set of natural numbers.

At the same time, the name of a set unites all elements of this set into one entity – this set. Consequently, all elements from the set $X$ have one common name "an element from $X$". This name discerns elements from $X$ and all other entities.

The second approach to set definition, intentional method gives rules for constructing elements that form some set. Consequently, each element in such set is named, in some sense, and its name is the procedure (algorithm, program) by which the element in question is constructed. This set is actually a single named set. In it, all elements that belong to the set with this name ($M$, for example) have the same, common name ("an element from the set $M$"). In such a way, we find that all ordinary sets are special kinds of named sets. Namely, they are singlenamed sets when we consider their structure in more detail.

In this setting, a mapping of ordinary sets is a morphism of the corresponding named sets. This allows one to understand and reflect more completely the structure of sets and relations between them. Such new expressive opportunities are especially important for applications. For example, this technique is used to solve some paradoxes in biological classifications (Burgin, 1983).

Let us consider in more detail how singlenamed sets such that come from the set theory emerge in the axiomatic set theory. According to Fraenkel and Bar-Hillel (1958), the most important directions in the axiomatic set theory are Zermelo-Fraenkel (**ZF**), von Neumann (**VN**), Bernays-Godel (**BG**), Quine (**NF**, **ML**) and Hao Wang (**Σ**) theories.

The **ML** theory of Quine is based on the concept of a class. There are two kinds of classes. However, the question if all objects in **ML** are classes or there are elements that are not classes.

In axiomatic set theory, any name is an expression of the corresponding set theory language, for example, a single sign (e.g., *M*) or some logical formula $\{x \in X; P(x) \,\&\, (y \in Y ((x, y) \in A \subseteq X \times Y)\}$. This gives a formal description of a set in a set theoretical language using set variables, logical signs and other mathematical symbols. As an example, we can take such axiomatic theories as Zermelo-Fraenkel (**ZF**) or Bernays-Godel (**BG**) theories.

However, set is a single primitive element only in some theories of sets. As an example we can take the axiomatic system **ZF** of Zermelo and Fraenkel or the universe of constructive sets built by Gödel (1940). Such theories as **NF** and **ML** work wlth classes instead of sets. In this case, instead of singlenamed sets, we have singlenamed classes.

Thus, we come to the conclusion that *any set is a named set and any class is a named class* because we cannot separate a name from the set or class in question. In addition, several other named sets (named classes) are connected with a set (class). They reflect different aspects of the initial set (class).

Another fundamental structure of the whole mathematics is **function**. However, functions are special kinds of binary relations between two sets, i.e., such relations in which any element from the first set is corresponded to at most one element from the second set.

**Definition 3.2.** *A function from X to Y is defined* (e.g., (Herrlich, and Strecker, 1973)) *to be a tripple* $(X, f, Y)$ *where* $f \subseteq X \times Y$ *is a binary relation such that for each* $x \in X$ *there is exactly one* $y \in Y$ *such that* $(x, y) \in f$.

This definition stems from the earlier approach to the definition of a correspondence, which was introduced by Bourbaki (1960):

**Definition 3.3.** *A correspondence between sets A and B is a triple* $(A, G, B)$ *where G is a graph, i.e., a set of pairs, for which* $\text{pr}_1 G \subseteq A$ *and* $\text{pr}_2 G \subseteq B$.

As Herrlich, and Strecker (1973, p. 10) mention, sometimes the triple $(X, f, Y)$ is denoted by $f: X \to Y$ or occasionally (and inaccurately) by $f$ alone. This explicitly shows that a function, as well as a correspondence, is a special case of named sets or fundamental triads. They are called set theoretical named sets, which are defined below.

**Definition 3.4.** *A named set* $\mathbf{X} = (X, r, I)$ *is called set theoretical if X and I are sets and r is a subset of the direct product* $X \times I$, *i.e.,* $r \subseteq X \times I$.

We see that a set theoretical named set coincides with a correspondence in the sense of Bourbaki. Some are confused with this and think that the concept of a named set is not new. However, there are named sets that are not set theoretical. Algorithms give examples of such named sets (cf. Section 3.5). Multisets (Aigner, 1979; Hickman, 1980; Knuth, 1997; 1998) and multigraphs (Berge, 1973) give other examples of named sets that are not set theoretical.

Thus, we conclude that **any set and any function are named sets**.

### 3.2. Generalization of Sets

In spite that sets are usually taken as the most fundamental entity in mathematics, there are different generalizations of sets. Let us consider such generalizations and their relations to named sets. We begin with multisets as they appeared much earlier than other generalizations of sets. Knuth (1998) attributes introduction of multisets to the Indian mathematician Bhascara Acharya (circa 1150), who studied permutations of multisets.

It is interesting to note that necessity in this structure is so urgent that multisets have been several times rediscovered and appeared in literature under different names (Blizard, 1991). For example, they are called bags by Peterson (1981), who attributed this term to Cerf, Gostelow, Estrin, and Volanski (1971).

**Definition 3.5** (Knuth, 1997). *A multiset is a collection that is like a set but can include identical or indistinguishable elements*.

For instance, $\{a,a,b,b,b\}$ is a multiset that contains two elements $a$ and three elements $b$. Thus a multiset is obtained if we take a named set **X** and add an axiom demanding that elements from the support S(**X**) of **X** are distinguishable if and only if they have different names in **X**.

Another interpretation of a multiset emphasizes that it allows repetitions of the same object.

**Definition 3.6** (Aigner, 1979). *A multiset on a set S is a function $r: S \to N$ that defines multiplicity of the elements from S* (*here $N = \{0, 1, 2, ... , \}$*).

Such multisets are named sets for the case when **Ens** consists of arbitrary sets and their maps, **Set** contains the single object $N$ and all binary relations on it, while **Col** is equal to **Ens**.

**Definition 3.7** (Hickman, 1980). *A multiset on a set S is a function $r: S \to$ **Card** that defines multiplicity of the elements from S*. Here **Card** is the class of all cardinals.

So **multisets**, in this extended sense, **are also special cases of named sets**.

Now let us consider fuzzy sets. To make this paper self-consistent, we give some definitions from fuzzy sets theory because some concepts have different interpretations. Fuzzy sets were introduced by Zadeh in 1965. The aim was to get better models for real-life systems and processes than allowed set theory. At the same time, in an abstract algebraic context, Salii (1965) defined a more general kind of structures, which were called L-relations.

**Definition 3.8** (Zadeh,1965). *A fuzzy set A in a set U is the triad* $(U, \mu_A, [0,1])$, *where $\mu_A: U \to [0,1]$ is a membership function of A and $\mu_A(x)$ is the degree of membership in A of $x \in U$.*

So, **fuzzy sets** in the sense of Zadeh **are the named sets** in the case when the category **Ens** consists of arbitrary sets and their mappings, while the category **Set** contains the single object [0,1] and all relations on it.

Different authors introduced more general kinds of fuzzy sets.

**Definition 3.9** (McVicar-Whelan, 1977). *A fuzzy set A in a set U is the triad (U, $\mu_A$, [-1,1]), where $\mu_A: U \to [-1,1]$ is a membership function of A and $\mu_A(x)$ is the degree of membership in A of $x \in U$.*

So, **fuzzy sets** in the sense of McVicar-Whelan **are the named sets** in the case when the category **Ens** consists of arbitrary sets and their maps, while the category **Set** contains the single object [-1,1] and all relations on it.

**Definition 3.10** (Cai Wen, 1984). *A fuzzy set A in a set U is the triad (U, $\mu_A$, (-$\infty$, +$\infty$) ), where $\mu_A: U \to (-\infty, +\infty)$ is a membership function of A and $\mu_A(x)$ is the degree of membership in A of $x \in U$.*

So, **fuzzy sets** in the sense of Cai Wen **are the named sets** in the case when the category **Ens** consists of arbitrary sets and their maps, while the category **Set** contains the single object (-$\infty$, +$\infty$) and all relations on it.

Let $L$ be a complete lattice, i.e., a partially ordered set with operations sup and inf.

**Definition 3.11** (Gogen, 1967). *An L-fuzzy set A in a set U is a triad (U, $\mu_A$, L), where $\mu_A: U \to L$ is a membership function of A and $\mu_A(x)$ is the degree of membership in A of $x \in U$.*

When $U \subseteq X \times Y$, we have *L*-fuzzy relation in the sense of Salii (1965).

So, ***L*-fuzzy sets** in the sense of Gogen **are the named sets** in the case when the category **Ens** consists of arbitrary sets and their maps, while the category **Set** contains the single object *L* and all relations on it.

**Definition 3.12** (Zimmermann, 1991). *A fuzzy set A in a set U is a triad (U, $\mu_A$, M), where M is the membership space, $\mu_A: U \to M$ is a membership function of A and $\mu_A(x)$ is the degree of membership in A of $x \in U$.*

**Remark 3.2.2.** The definition 3.2.8, in which *M* may be an arbitrary set, does not preserve the main idea of fuzzy set theory that the values of the membership function $\mu_A(x)$ reflect to what extent elements from the universal set *U*, which is the domain of $\mu_A(x)$, belong to the fuzzy set A. To reflect the extent of membership, the codomain *I* of $\mu_A(x)$ has to be, at least, partially ordered. At the same time, structures that have the form A = (*U*, $\mu_A$, *I* ), and in which *I* is an arbitrary set are considered in the theory of named sets and are called set-theoretical named sets (Section 3.1).

Intuitionistic fuzzy sets, which were introduced by Atanasov [3, 4], are represented either by two named sets or by a named set in which the naming relation is not a function.

In more detail, relation between fuzzy sets and named sets are considered in (Burgin and Kuznetsov, 1992). In any case, we see that **all kinds of fuzzy sets and fuzzy relations are special cases of named sets**.

It is also demonstrated (Burgin and Chunihin, 1997; Chunihin, 1997) that named sets provide more extended means than fuzzy sets for representing and investigating uncertainty.

### 3.3. Mathematical logic

The subject of mathematical logic has origins in philosophy, in general, and in philosophy of mathematics, in particular. Indeed, it is essentially nonmathematical argument that can show the usual rules for inference and/or deduction (such as the law of excluded middle; cut rule; etc.) are valid. It is also a legacy from philosophy that we can distinguish semantic reasoning (answering the question "what is true?") from syntactic reasoning (answering the question "what can be shown to be true?"). The first question leads to model theory, the second one, to proof theory. In spite of this, mathematical logic is often considered as the study of the processes used in mathematical deduction, ignoring model theory. Actually, we have two traditional parts of mathematical logic: syntax (proof theory) and semantics (model theory). Pragmatics of mathematical logic is only emerging.

Typically, a logic **L, as a mathematical structure,** consists of a formal or informal language *L* together with a deductive system and/or a model-theoretic semantics. The language is, or corresponds to, a part of a natural language like English or Greek. The deductive system is developed to capture, codify, or simply record which *inferences* are correct for the given language, and the semantics is built to capture, codify, or record the meanings, in the form of truth-conditions, or possible truth conditions, for at least part of the language. We call it a logical semantics.

A deduction system consists of deduction rules. Each deduction rule and consequently, each step of deduction has the form

$$A \to B \qquad (13)$$

Here *A* and *B* are some expressions or sets of expressions from the language of the corresponding logic.

Thus, a **deduction rule is the named set** $(A, \to, B)$.

Likewise, any formal calculus **C** is also a special kind of named sets. Really, it may be considered as a triad **C** = (*A, R, T*) where *A* is the set of axioms, *R* are rules of deduction by which from axioms the theorems of the calculus are deduced. These theorems form the set *T*. The named set (*A, R, T*) is called a named set of calculus rules. The same calculus may be represented by another (deduction) named set (*A, d, T*) where the relation *d* connects any axiom *a* from *A* with such theorems *t* from *T* that *a* is used in a process of deduction of *t*.

A logical semantics corresponds truth-values to the sentences of the language *L*. That is, if we consider a language as a set of its sentences, then the semantics of a logic **L** is partial function from *L* into the set TV of the truth-values of the logic **L**. For example, in classical logics TV consists of two symbols {1, 0}, while in fuzzy logics (Zimmermann, 1991) TV is the unit interval [1, 0]. As it has been demonstrated in section 3.1, any partial function is a named set. Consequently, **any logical semantics is a named set**.

It is possible to show that any language, formal or natural, is built of different named sets. All this demonstrates that **any logic is built of different named sets**.

In addition, the same is true for such an important logical construction as model. To explain this, it should be noted that there exists no generally accepted and exact definition of a model in mathematical logic. Some authors treat a model simply as a mathematical structure, exactly as a set with a system of relations (may be functions) in it (Malzev, 1970; Shoenfield, 1967). Other authors interpret a model as a pair consisting of a mathematical structure (the same as above) and a partial mapping of a logical language into this structure (Chang and Keisler, 1973; Mendelson, 1963). In the second case, we again have some named set because any partial mapping is a named set. Really, speaking about a model of some logical language, one factually bears in mind the named set ($L, i, M$) of an interpretation of the language $L$ into some mathematical structure $M$. A mathematical structure $M$ (that is called a model of the language $L$) is the set on which the predicates and functions are defined. The mapping $i$ is built in such a manner that there are three correspondences:

(i) between predicate symbols from the alphabet of the language $L$ and predicates having the same number of variables defined in $M$;

(ii) between functional symbols and functions on $M$;

(iii) between constants and elements of the mathematical structure $M$.

Thus, **any model is a named set**.

In a similar way, an interpretation of a formal theory (or a deductive calculus) $T$ into a model $M$ is built. The corresponding modeling named set has a form ($T, p, M$) and the truth is such a property that its conservation is demanded from the naming relation $p$. In other words, if all formulae from $T$ are considered as true, then $M$ is a model of $T$ if and only if the images of these formulae are true in $M$. It should be noted that when in logic (considered as a model), it is taken a pair that consists of the interpretation map and a mathematical structure, then such model is the part of the modeling named set.

### 3.4. Theory of Categories

Category theory is a highly abstract mathematical theory of structures and, even more, of systems of structures. It allows us to study, among other things, how structures of some kind as well as systems of structures of different kinds are related to one another. As a result, category theory unifies and provides a fruitful organization of mathematics. Results that were obtained for categories are applied to diverse mathematical structures such as groups, sets, topological spaces, etc. In particular, category theory and even its part, which is called topos theory (Goldblatt, 1979), is considered as a new foundation for mathematics. The only problem is, as states Stanford Encyclopedia of Philosophy, whether it should be considered seriously as providing a foundational alternative to set theory or whether it is foundational in a different sense altogether.

According to the most popular approach, we have the following definition of a category.

**Definition 3.13.** *A category* **C** *consists of two collections* **Ob C**, *the objects of* **C**, *and* **Mor C**, *the morphisms of* **C**, *which satisfy the following three axioms*:

- **A1.** For every pair $A$, $B$ of objects, there is a collection $\mathbf{Mor}_C(A, B)$, which is a part of **Mor C.** Elements from $\mathbf{Mor}_C(A, B)$ are called morphisms from $A$ to $B$ in **C**. When $f$ is a morphism from $A$ to $B$, it is denoted by $f: A \to B$.
- **A2.** For every three objects $A$, $B$ and $C$ from **Ob C**, there is a partial operation, which is a partial function from pairs of morphisms that belong to the direct product $\mathbf{Mor}_C(A, B) \times \mathbf{Mor}_C(B, C)$ to morphisms in $\mathbf{Mor}_C(A, C)$. In other words, when $f: A \to B$ and $g: B \to C$, there is a morphism $g \circ f: A \to C$ is called the composition of morphisms $g$ and $f$ in **C**. This composition is associative, that is, if $f: A \to B$, $g: B \to C$ and $h: c \to D$, then $h \circ (g \circ f) = (h \circ g) \circ f$.
- **A3.** For every object $A$, there is a morphism $1_A$ in $\mathbf{Mor}_C(A, A)$, called the identity on $A$, for which if $f: A \to B$, then $(1_B \circ f) = f$ and $(f \circ 1_A) = f$.

Examples of categories:

1. The category of sets **SET**: in it objects are arbitrary sets and morphisms are mappings of these sets.
2. The category of groups **GRP**: in it objects are arbitrary groups and morphisms are homomorphisms of these groups.
3. The category of topological spaces **TOP**: in it objects are arbitrary topological spaces and morphisms are continuous mappings of these topological spaces.

In another approach, a category **C** consists only of the collection **Mor C** of the morphisms of **C**. Objects of **C** are associated with identity morphisms $1_A$. It is possible to do because $1_A$ is unique in each set $\mathbf{Mor}_C(A, A)$, and uniquely identifies the object $A$. In any case, morphism is the central concept in a category. But what is a morphism?

If $f: A \to B$ is a morphism, then it is a one-to-one named set $(\{A\}, f, \{B\})$. Thus, the main object of a category is a named set and categories are build from these named sets. Besides, a construction or separation of a category begins with separation of all elements into two sets and calling all elements from one of these sets by the name "objects" and elements from the other sets by the name "morphisms". In such a way, two named sets appear. In addition, composition of morphisms, as any algebraic operation, is also represented by a named set ( $\mathbf{Mor}_C(A, B) \times \mathbf{Mor}_C(B, C)$, $\circ$ , $\mathbf{Mor}_C(A, C)$ ). This shows that even if categories and sets are used for formalization of the concept of a named set, the informal notion of a named set is prior both to categories and sets.

As a result, we come to the conclusion that **any category is built of different named sets**.

Moreover, functors between categories, which are structured mappings of categories (Herrlich and Strecker, 1973), are morphisms of those named sets. Thus, not only categories are used for formal definition of named sets, but also informal named sets are used for definition of categories and relations between them.

In more detail, the correspondence between categories and named sets is studied in (Burgin, 1991; 1991a).

### 3.5. Natural Numbers

Natural numbers are often considered as the most fundamental objects in mathematics. As wrote such prominent mathematician as Leopold Kronecker, "*God made the integers, all the rest is the work of man*".

There are several approaches to introduction of natural numbers in mathematics.
1. As it is done in counting and then formalized in the Peano axioms for arithmetic, natural numbers are constructed by taking the next element. This step of going from $n$ to $n + 1$ or to $Sn$ is represented by a named set $(n \to n + 1)$ or $(n \to Sn)$.

2. When mathematics is based on set theory, natural numbers are defined and treated as equivalent classes of finite sets. Thus, we natural come to the named set in which the support consists of all finite sets, as its elements, and the set of names comprises denotations of natural numbers. For example, if the decimal system is used, then all sets with 10 elements are corresponded to the number 10, while if the binary system is used, the same sets are related to the number 1010.
3. Another approach to natural numbers considers them as properties of sets (Burgin, 1989b). This approach is formalized in the theory of abstract properties (Burgin, 1989).

**Definition 3.14.** *An abstract property* **P** *has the form of the named set* **P** = (*U*, *p*, *L*) *where U is a universe where the property* **P** *is defined, L is the scale of the property* **P**, *and p*: $U \to L$ *is a partial function that corresponds to objects from U elements from L.*

Thus, we encounter a new named set:

$$U \xrightarrow{p} L \qquad (14)$$

For example, the most popular property in logic is truth. For this property, the universe $U$ consists of statements or propositions of a logical language, $L = \{T, F\}$ in classical logics, and $p$ corresponds to sentences T when these sentences are true and F when these sentences are false. Another example is from physics. For such property as mass, the universe $U$ consists of all material bodies, $L$ is equal to the set $\mathbf{R}^+$ of non-negative real numbers, and $p$ corresponds to each material body its mass.

To consider natural numbers as an abstract property, we take the collection of all sets that exist in reality as the universe on which the property is defined. We denote this universe by $W$. The scale $L$, in this case, is some denotational system for natural numbers. For example, $L$ may be the decimal system or binary representations of natural numbers. In any case, $L$ is a set of words in some alphabet $A$ and these words are names of natural numbers. In our examples, $A$ is either $\{0, 1, 2, 3, \ldots, 9\}$ or $\{0, 1\}$. In the first case, we have such names as 123, 54, 1111 and so on. In the second case, we have such names as 101, 10, 1111 and so on. For distinction, we denote such a scale for natural numbers by $N$. To a set $X$ from $W$, the mapping $p$ corresponds the name of the number of the elements in *X*. This is exactly the way in which natural numbers appear in reality. For example, *p*({a, b, c}) = 3 and *p*({b, b}) = 2. For distinction, we denote this mapping by η. In such a way, we obtain the abstract property ***N*** = (*W*, η, *N*). We call this property a natural number property. Depending on the scale *N*, we have as many natural number properties as there are natural number systems. For example, when *N* consists of all words in the alphabet {0, 1, 2, 3, … , 9 }, we have the decimal number system. When *N* consists of all words in the alphabet {0, 1}, we have the binary number system, which is used in computers.

In this context, a separate number is not a property, but a value of the natural number property, while all numbers are represented as a property. When we say that, for example, 1 is a natural number, it is not exact because it is only a name of a natural number. More exactly, we have to speak about a natural number with the name 1. However, in natural languages, it is correct to speak about 1 as a natural number as meaning of this statement is a natural number with the name 1.

Thus, **natural numbers** form a property of sets, which in its turn, **is a named set**.

This model of natural numbers gives a supporting argument for the philosophical approach of mathematical realism because natural numbers (or finite cardinals) exist exactly in this form.

### 3.6. Theory of Algorithms

Algorithms are used as used as a basic object for foundations of mathematics in constructive approaches in mathematics (Markov, 1962; Bishop and Bridges, 1985; Beeson, 1985). Consequently, theory of algorithms forms one more basis for mathematics. In addition, theory of algorithms is one of the foundations for computer science (Aho and Ullman, 1994). Let us analyze its main concepts.

A popular point of view on algorithm is presented by Rogers (1987):

Algorithm is a *clerical* (*i.e., deterministic, book-keeping*) procedure which can be applied to any of a certain class of symbolic *inputs* and which will eventually yield, for each such input, a corresponding *output*.

More generally, an algorithm is a specific kind of recipe, method, or technique for doing something. In other, more exact words, *an algorithm is a text giving unambiguous*

*(definite) and simple to follow (effective) prescriptions (instructions) how from given inputs (initial conditions) derive necessary results*.

However, this notion of algorithms is imprecise and, consequently, insufficient for a mathematical study of algorithms. Consequently, mathematical models of algorithm have been elaborated and corresponded to informal notions of algorithm. In contrast to real life algorithms, which are applied to a multiplicity of different kinds of objects, mathematical algorithms work only with symbolic expressions, which are words of some formal language.

The first step in constructing a formal model is specification of an algorithm, which is developed as a formal counterpart to the idea of algorithm:

- *A uniform specification of the symbolic expressions accepted as sets of instructions and of those objects* (*expressions, things, actions, processes, etc.*) *that are considered as inputs, and as outputs* (EP-formalism).
- *A uniform specification how instructions and input determine subsequent functioning of algorithms*.
- *A uniform specification how the output of a given computation is obtained and determined*.

These specifications show that any algorithm $A$, in addition to instructions or rules $R$, contains the input language $L_I$ and the output language $L_O$. Thus, an algorithms is related to a named set $A = (L_I, R, L_O)$. These named sets are not set-theoretical because different sets of rules may determine the same function from $L_I$ to $L_O$.

All three components of the algorithm's named set $A$ in their turn consist of other named sets. A named set that corresponds to an arbitrary formal language is described in (Burgin, 1997). For rules $R$ of the algorithm, this is shown below in the context of formal models of algorithm.

Aiming at solution of different problems, researchers have suggested a diversity of exact mathematical models for a general notion of algorithm. Recursive functions (Gödel, 1934), ordinary Turing machines (Turing, 1936; Post, 1936), λ-calculus (Church, 1932/33; 1936), recursive and partial recursive functions (Kleene, 1936) appeared the first.

Then other mathematical models of algorithms were suggested. They include a variety of Turing machines appeared: Turing machines with several tapes, with several heads, with n-dimensional tapes, non-deterministic, probabilistic, alternating Turing machines, etc. The most popular of other mathematical models are: neural networks; finite automata; formal grammars; Petri nets; Post productions; Minsky machines; random access (RAM), RASP (Random Access Stored Program), vector, array and hardware modification machines.

Here we consider only some of these models. We begin with one of the most popular models, Turing machine. A Turing machine **T** is a triad

$$T = (L, D, Q)$$

where **L** is the *language* of **T**, **D** is the *device* of **T**, and **Q** is the *state configuration* of **T**.

The *language* **L** of the Turing machine **T:**

$L = ( L_I , L_W , L_O )$ where $L_I$ is the *input* language, $L_W$ is the *working* language, and $L_O$ is the *output* language of **T**. $L_X = (A_X, R_X, L_X)$ where $A_X$ is the *alphabet*, $R_X$ is the set of *generating rules*, and $L_X$ is the *set of all words* of $L_X$. Usually $L_X$ is the set of all finite strings in the alphabet $A_X$.

The *device* **D** of the Turing machine **T:**

$D = ( H, P, M)$ where P is the *program*, H is the *operating device* (which is called traditionally the head of **T**), and M is the *memory* of **T** (which consists traditionally of one or several tapes, which may be expanded indefinitely, involving, thus, a potential infinity).

The *state configuration* **Q** of the Turing machine **T:**

$Q = (q_0 , Q , F)$ where Q is the *state space*, $q_0$ is the *initial state*, and F is the *set of final states* of **T.**

The *program* P of the Turing machine **T** consists of rules that prescribe actions of the Turing machine and have one of the following forms:

$$q_h r_i \to r_j q_k, \quad q_h r_i \to R q_k \quad, \quad q_h r_i \to L q_k$$

In deterministic Turing machines, there are no two rules in P with the same left part. In nondeterministic Turing machines, one and the same state and symbol may lead to different actions of a Turing machine.

Thus, we come to a conclusion that a Turing machine is structured as a hierarchy of triads and the central element of a Turing machine, its rule is a fundamental triad/named sets

Another popular model of algorithms is finite automaton.

A finite automaton **A** consists of three structures:

- The *linguistic structure* L = ( $\Sigma$, Q, $\Omega$ ) where $\Sigma$ is a finite set of *input symbols*, Q is a finite *set of states*, and $\Omega$ is a finite set of *output symbols* of the automaton **A**;

- The *state structure* S = ( Q, $q_0$ , F ) where $q_0$ is an element from Q that is called the *start state* and F is a subset of Q that is called the set of *final* (in some cases, *accepting*) states of the automaton **A**;

- The *action structure*, which is traditionally called a *transition function*

$$\delta: \Sigma \times Q \to Q \times \Omega$$

Thus, we can see **A** is the triad (L, S, $\delta$) and all its components, L, S, and $\delta$, are named sets or fundamental triads.

One more model of algorithm is generative grammar. In the most general form, it is called unrestricted grammar, or phrase-structure grammar.

An *unrestricted grammar*, or *phrase-structure grammar*, **G** has the form

$$( L, S, P )$$

where

- L = ( V, $\Sigma$ ) is the *lexical structure* of **G** where V is a finite *set of variables* and $\Sigma$ is a finite set of *terminal symbols*, or *terminals*;

- S is the *start symbol* of **G**;

- P is the *generative structure* of **G**, which consists *productions* of the form

$$\alpha \to \beta$$

where $\alpha$ and $\beta$ belong to the set ( V$\cup\Sigma$ )* of all strings in V$\cup\Sigma$ and $\alpha$ contains at least one variable.

Thus, we can see that **G** is the triad (L, S, P) and all its rules $\alpha \to \beta$ are named sets or fundamental triads.

Partially recursive functions are functions and thus, named sets (cf. Section 3.1).

Post productions have the same structure as rules of generative grammars and thus, they are named sets.

All formulas are built of symbols. So, concept of symbol is basic for the constructive and logical foundations of mathematics. At the same time, it is demonstrated in (Burgin, 1993a) that different models of symbols are built of different named sets.

## 4. Conclusion

Thus, we have demonstrated that named set is the basic structure in mathematics providing for the most fundamental foundations of mathematics and allowing systematizing and analyzing all other foundations.

In addition, discovery of named sets made possible to explicate and explain main directions in the mathematics development. The progress of contemporary mathematics may be viewed as a process of transition from ordinary (singlenamed) sets to more general cases. To show this, we can consider the main mathematical fields: mathematical logic, theory of algorithms, algebra, geometry, topology, arithmetic and number theory, combinatorics, functional analysis, theory of differential equations, and probability.

For example, in topology a lot of achievements is connected with the introduction of fibers and their special cases fiber spaces, bundles, smooth fibers etc. But any fiber **F** is a topological functional named set **F** = (*E, p, B*), i.e., a named set in which the base *B* and the fiber space *E* are topological spaces and *p* is a continuous projection of *E* onto *B*.

In many cases, fibers, which are special cases of topological named sets, replace topological spaces in topology. Manifolds are used instead of Euclidean spaces in mathematical analysis. If we take modern analysis the main structure on which functions are defined and studied are different kinds of manifolds: topological, differentiable, smooth etc. The most general kind is a topological manifold. Each of these structures is a named set **X** = (*X, r, I*) where *X* is a topological space, *I* is some n-dimensional Euclidean space $E^n$ and *r* is a continuous relation that is a local homeomorphism, i.e. for each point *x* from *X* there exists an open neighborhood that is homeomorphic to some open subset of $E^n$. Conditions that define special cases of topological manifolds (differentiable, analytical etc.) may be also formulated in the language of the named set theory, i.e., as conditions on named sets that are obtained by application to these named sets definite set-theoretical operations (cf. Burgin, 1984).

Scalar, vector and tensor *fields on manifolds are also named sets* having the form **X** = (*X, n, D*) where *X* is the same as above, *D* is some set of scalars (real or complex numbers), vectors or tensors and *n* is a function defined on *X*.

It is necessary to remark that transition from a set to a named set as a basic structure (a carrier for other structures) for the development of different fields in mathematics is peculiar to the majority of mathematical fields. In algebra, in addition to ordinary (or homogeneous) universal algebras, in which operations are defined on a set, heterogeneous algebras, in which operations are defined on a named set, are studied.

**Definition 4.1**. *A heterogeneous or multibase universal algebra* **A** *is a set A with a system of operations* Σ *in which elements of A form an indexed system A = {$A_i$} of sets and each operation is a mapping having the form f: $A_{i1} \times A_{i2} \times \ldots \times A_{ik} \to A_i$*.

**Definition 4.2**. *The system A is called the carrier of the multibase algebra* **A**.

Examples of multibase universal algebras are modules, polygons, i.e., sets on which monoids act, automata, polyadic or Halmos algebras (Halmos, 1962), nonhomogeniuos polyadic algebras (Leblanc, 1962), relational algebras (Beniaminov, 1979) and state machines. Multibase universal algebras were studied by several authors under different names. To mention only some of them, it is necessary to name algebras with a scheme of operators introduced by P.J.Higgins (1963; 1973), heterogeneous algebras from the papers of Birkhoff and Lipson (1970) and Mathienssen (1978), and many-sorted algebras studied by Plotkin (1991). The term " heterogeneous algebras" is used more often than other related terms. Heterogeneous (multibase or many-sorted) algebras represent the next level of the development of algebra. Namely, in ordinary (or homogeneous) universal algebras operations are defined on a set, while in heterogeneous algebras operations are defined on a named set. This makes possible to develop more adequate models for many processes and systems. For example, heterogeneous algebras are extensively used for mathematical modeling information processing by computers. Such models as abstract automata and abstract states machines or evolving algebras become more and more widespread in computer science (Gurevich, 1995). In addition, relational algebras are extensively used for modeling relational databases (Beniaminov, 1979; 1980; Plotkin, 1991). Modern combinatorics is built not on sets but on multisets, which are special cases of named sets. Multivalued and multi-sorted logics, semantics of which is built not on sets but on named sets, are becoming more and more popular in logic.

Demonstration that many mathematical structures are either named sets or built of named sets supports the main thesis of this work that named set is the most fundamental structure in mathematics basic for its foundations.